\documentclass[11pt]{amsart}
\usepackage{amsmath,amsthm,amssymb,latexsym,color}
\usepackage{hyperref}
\usepackage{lipsum}

\theoremstyle{plain}
\newtheorem{theor}{Theorem}[section]

\newtheorem{prop}[theor]{Proposition}

\theoremstyle{remark}

\renewcommand{\a}{\alpha}
\renewcommand{\b}{\beta}

\def\PP{{\mathbb P}}
\def\E{{\mathbb E}}

\def\R{{\mathbb R}}

\def\N{{\mathbb N}}
\def\Z{{\mathbb Z}}
\def\Event{{\mathcal E}}

\def\Prob{{\mathbb P}}

\newcommand{\brk}[1]{\left( #1\right)}
\newcommand{\brkb}[1]{\big( #1\big)}
\newcommand{\set}[1]{\left\{#1\right\}}
\newcommand{\abs}[1]{\left| #1 \right|}
\newcommand{\norm}[1]{\left\Vert #1 \right\Vert}
\newcommand{\normb}[1]{\big\Vert #1 \big\Vert}

\title{A remark on the smallest singular value of powers of Gaussian matrices}

\author{Han Huang}
\author{Konstantin Tikhomirov}

\address{School of Mathematics, Georgia Institute of Technology}
\email{hhuang421@gatech.edu}
\email{konstantin.tikhomirov@math.gatech.edu}

\begin{document}

\begin{abstract}
Let $n,k\geq 1$ and let $G$ be the $n\times n$ random matrix with i.i.d.\ standard real Gaussian entries.
We show that there are constants $c_k,C_k>0$ depending only on $k$ such that the smallest singular
value of $G^k$ satisfies
$$
c_k\,t\leq \Prob\big\{s_{\min}(G^k)\leq t^k\,n^{-1/2}\big\}\leq C_k\,t,\quad t\in(0,1],
$$
and, furthermore,
$$
c_k/t\leq \Prob\big\{\|G^{-k}\|_{HS}\geq t^k\,n^{1/2}\big\}\leq C_k/t,\quad t\in[1,\infty),
$$
where $\|\cdot\|_{HS}$ denotes the Hilbert--Schmidt norm.
\end{abstract}

\maketitle

\section{Introduction}

Everywhere in the paper, $G$ denotes an $n\times n$ random matrix with i.i.d.\ real valued standard Gaussian entries. 
The smallest singular value and the condition number of standard square Gaussian matrices (and other random matrix models)
are classical objects of interest within the random matrix theory.
The condition number $\kappa(A)=s_{\max}(A)/s_{\min}(A)$ of a matrix $A$ is of importance as a simple estimator of the relative error
when solving the linear system $Ax=b$ with the coefficient vector $b$ known up to some additive error
(see, for example, \cite{Smale}).

In 1940-es, von Neumann and
Goldstine \cite{vNG} conjectured that the ``typical'' value of $s_{\min}(G)$
is of order $n^{-1/2}$, while the condition number $\kappa(G)=s_{\max}(G)/s_{\min}(G)$ is of order $n$.
The conjecture
was rigorously established by Edelman \cite{Edelman} and, independently, by Szarek \cite{Szarek}.
The proofs in \cite{Edelman,Szarek} use as the central element a formula for the joint distribution of singular values of $G$.
In particular, the following estimate for the smallest singular value of $G$ was obtained in \cite{Edelman,Szarek}:
\begin{equation}\label{eq:Szarek}
\Prob\big\{s_{\min}(G)\leq t\,n^{-1/2}\big\}=\Theta(t),\quad t\in(0,1].
\end{equation}
Here, we adopt the ``big theta'' notation:
given two non-negative functions $f(t)$ and $g(t)$ defined on the same domain, we write $f(t)=\Theta(g(t))$
if $C^{-1} f(t)\leq g(t)\leq Cg(t)$ for all $t$ and some universal constant $C\geq 1$.
When the constant is allowed to depend on a parameter, we add the parameter as a subscript for $\Theta$.
Numerous results dealing with invertibility of non-Gaussian random models have appeared in literature.
We prefer to avoid discussion of that (very active) research direction in this note.
Let us refer to surveys \cite{RV10} and \cite{Rud14} which give some (partial) account of the subject.

Returning to linear systems with random coefficients,
it seems natural to consider the situation when we are given a linear system of the form $G^k x=b$,
where $k\geq 1$ is fixed, and would like to estimate the relative error of the obtained solution when $b$ is known
up to some additive error.
In this case, we could ask what is the typical value of the condition number of $G^k$ and, moreover, what
are optimal large deviation estimates for $\kappa(G^k)$?
Since the largest singular value of $G^k$ is of order $\Theta_k(n^{k/2})$ with a very large probability,
the question essentially amounts to computing small ball probabilities for $s_{\min}(G^k)$.
Obviously, the trivial relation $s_{\min}(G^k)\geq (s_{\min}(G))^k$ and the known estimates for $s_{\min}(G)$
immediately imply probabilistic estimates for $s_{\min}(G^k)$, which, however, turn out to be suboptimal.
In this note, we are interested in non-asymptotic estimates which are sharp up to multiplicative constants.
To authors' best knowledge no such results have been previously noted in the literature.
The main statement of the note is
\begin{theor} \label{thm: main}
Let $n,k\geq 1$ and let $G$ be the $n\times n$ matrix with i.i.d.\ standard Gaussian entries. 
Then
\begin{align*}
&\Prob\big\{\|G^{-k}\|_{HS}\geq t^k\,\sqrt{n}\big\}=\Theta_k(1/t),\quad t\in[1,\infty),\quad\mbox{and}\\
&\Prob\big\{s_{\min}(G^{k})\leq t^k\,n^{-1/2}\big\}=\Theta_k(t),\quad t\in(0,1].
\end{align*}
\end{theor}
Here, $\|\cdot\|_{HS}$ denotes the Hilbert--Schmidt norm of a matrix.

\bigskip

{\bf Acknowledgement.} The second named author is partially supported by the Sloan Research Fellowship.
Both authors are grateful to Mark Rudelson for interesting discussions.

\section{Proof}

Our proof relies on the following simple observation.
Let
$$
G = U\Sigma V^\top
$$
be the singular value decomposition of $G$, so that $\Sigma$ is the (random) diagonal matrix
with the singular values of $G$ arranged in the non-increasing order on the main diagonal,
and $U,V$ are (random) orthogonal matrices.
Further, let $W$ be an $n\times n$ random orthogonal 
matrix uniformly distributed on ${\rm O}_n(\R)$ (with respect to the Haar measure),
which is independent from $\{U,\Sigma,V\}$.
Then, in view of the invariance of the Gaussian distribution under orthogonal transformations,
the matrix $WG$ is equidistributed with $G$, whence
\begin{align*}
\normb{G^{-k}}_{HS}\stackrel{d}{=}
    \normb{(WG)^{-k}}_{HS} &= \normb{ \brkb{V\Sigma^{-1}U^\top W^\top}^k}_{HS}\\
     &= \normb{ V \brkb{\Sigma^{-1}U^\top W^\top V}^{k-1}\Sigma^{-1}U^\top W^\top}_{HS}\\
     &= \normb{  \brkb{\Sigma^{-1}U^\top W^\top V}^{k-1}\Sigma^{-1} }_{HS}\\
     &=\normb{  \brkb{\Sigma^{-1}Q}^{k-1}\Sigma^{-1} }_{HS},
\end{align*}
where the random orthogonal matrix $Q:=U^\top W^\top V$ is uniformly distributed on ${\rm O}_n(\R)$ and is independent from $\Sigma,G$.
Similarly, we have
$$
s_{\min}^{-1}(G^k)=\|G^{-k}\|\stackrel{d}{=}
    \normb{(WG)^{-k}}=\normb{\brkb{\Sigma^{-1}Q}^{k-1}\Sigma^{-1}},
$$
where $\|\cdot\|$ denotes the spectral norm.
Thus, the problem of estimating the right tail of the distribution of $\|G^{-k}\|_{HS}$ (and
of $\|G^{-k}\|$)
can be viewed as a particular case of a more general question of studying the distribution of the matrix product
$\brk{TW}^{k-1}T$, where $T$ is a fixed diagonal matrix and $W$ is uniformly distributed on ${\rm O}_n(\R)$.
\begin{prop} \label{prop:HSbound}
    Let $T={\rm diag}(\tau_1,\dots,\tau_n)$ be an $n\times n$ fixed diagonal matrix with non-negative entries, and let $W$ be 
    a uniform random orthogonal matrix. Take any $k\in\N$.
    Then
\begin{itemize}

\item For any even positive integer $m$ and any $i,j\in[n]$ we have
\begin{align*}
        \E  \brk{\brkb{\brk{T W}^k T}_{ij}}^m
        &\le C_{k,m}\tau_i^m\tau_j^m  \sum_{\b \in [n]^{m(k-1)/2}}
        \brk{\prod_{\ell=1}^{m(k-1)/2}\tau_{\b_\ell}^2} n^{-km/2}\\
        &= C_{k,m}\tau_i^m\tau_j^m \|T\|_{HS}^{m(k-1)}n^{-km/2},
\end{align*}
where $C_{k,m}>0$ depends only on $k$ and $m$.

\item The expectation of the squared Hilbert--Schmidt norm of $\brk{TW}^{k}T$ satisfies
    $$
        \E_W\normb{ \brkb{TW}^{k}T}_{HS}^2 \le C_k n^{-k}\norm{T}_{HS}^{2(k+1)},
    $$
    where $C_k>0$ only depends on $k$.

\item For any $i\leq n$, denoting by $T(i,s)$ the diagonal matrix with the $i$--th diagonal entry equal to $s$
    and all other entries equal to the corresponding entries of $T$, we have
$$
\big|\big\{s\in [\tau_i/2,\tau_i]:\; \E  \brkb{\brkb{\brk{T(i,s) W}^k T(i,s)}_{ii}}^2\geq 
c_{k}n^{-k}\tau_i^{2k+2}\big\}\big|\geq \tau_i/4,
$$   
where $c_k>0$ may only depend on $k$.
\end{itemize}
\end{prop}
Let us postpone the proof of the proposition till the end of the section, and complete the proof of the main result of the paper.

\subsection{Proof of Theorem~\ref{thm: main}}

With the matrices $\Sigma$ and $Q$ defined as above, application of Proposition~\ref{prop:HSbound} with $T:=\Sigma^{-1}$ and with $Q$
in place of $W$ gives 
\begin{align}
    \E_W\normb{(WG)^{-k}}^2_{HS}=
    \E_Q\normb{ \brkb{\Sigma^{-1}Q}^{k-1}\Sigma^{-1} }_{HS}^2 
    \le C_kn^{1-k} \normb{\Sigma^{-1}}_{HS}^{2k}= C_kn^{1-k} \normb{G^{-1}}_{HS}^{2k}.
    \label{eq: HSComparison}
\end{align}
It is clear that non-asymptotic estimates for the Hilbert--Schmidt norm of the inverse of the standard Gaussian matrix 
can be obtained by analysis of the joint distribution of its singular values, similar to \cite{Edelman, Szarek}.
However, we were not able to locate a ``ready-to-reference''
result of this kind in the literature, and instead will use
a more general statement about the Hilbert--Schmidt norm of the inverse of a random matrix with i.i.d.\ entries
with a continuous distribution \cite[Theorem~1.1]{Tikhomirov}, which implies, in particular, that
\begin{equation}\label{thm: G_HS}
    \PP\big\{ \norm{G^{-1}}_{HS}\ge tn^{1/2}\big\} \le \frac{C_{\text{\tiny{\ref{thm: G_HS}}}}}{t},\quad t>0,
\end{equation}
for a universal constant $C_{\text{\tiny{\ref{thm: G_HS}}}}\geq 1$.
Now, using \eqref{eq: HSComparison} and \eqref{thm: G_HS}, it is easy to
obtain the required upper bound on the right tail of $\|G^{-k}\|_{HS}$.

\smallskip

For $t>0$ and $i\in \Z$, let $\Event_i(t)$ be the event that $2^{i}tn^{1/2}\le \norm{G^{-1}}_{HS}\le 2^{i+1}tn^{1/2}$. Clearly,
\begin{align*}
    \PP \big\{ \normb{(WG)^{-k}}_{HS}\ge t^kn^{1/2}\big\}
    & =    \sum_{i\in \Z} \PP\big\{ 
        \normb{(WG)^{-k}}_{HS}\ge t^kn^{1/2} \,|\, \Event_i(t)\big\}\,
        \PP\brk{ \Event_i(t)}.
\end{align*}
For $i\ge 0$, by~\eqref{thm: G_HS} we have $\PP\brk{\Event_i(t)} \le \frac{C_{\text{\tiny{\ref{thm: G_HS}}}}}{2^it}$
and thus 
\begin{align}
\sum_{i = 0}^\infty  \PP\big\{ 
    \normb{(WG)^{-k}}_{HS}\ge t^kn^{1/2} \,|\, \Event_i(t)\big\}\,
    \PP\brk{ \Event_i(t)} \le \frac{2C_{\text{\tiny{\ref{thm: G_HS}}}}}{t}.
    \label{eq: tailPart}
\end{align}
For $i< 0$, everywhere on $\Event_{i}(t)$ we have $\norm{G^{-1}}_{HS}^{2k}\le 2^{2k(i+1)}t^{2k}n^{k}$.
Hence, conditioning on $\Event_{i}(t)$
and applying Markov's inequality together with \eqref{eq: HSComparison}, we obtain
\begin{align*}
    \PP_W\big\{\normb{(WG)^{-k}}_{HS}\ge t^kn^{1/2}\,|\,\Event_{i}(t)\big\}
    \le& t^{-2k}n^{-1}C_kn^{1-k} \cdot 2^{2k(i+1)}t^{2k}n^{k}= C_k2^{2k(i+1)},
\end{align*}
whence, again applying~\eqref{thm: G_HS},
\begin{align}
    \sum_{i=-\infty}^{-1} \PP\big\{\normb{(WG)^{-k}}_{HS}\ge t^kn^{1/2} \,|\, \Event_{i}(t)\big\}\,\PP\brk{\Event_{i}(t)}
    \le \sum_{i=-\infty}^{-1} C_k 2^{2k(i+1)}\cdot\frac{C_{\text{\tiny{\ref{thm: G_HS}}}}}{2^i t} 
\le \frac{C_{\text{\tiny{\ref{thm: G_HS}}}}\,C_k\,2^{2k}}{t}.
    \label{eq: MarkovPart}
\end{align}
Combining \eqref{eq: tailPart} and \eqref{eq: MarkovPart}, we obtain
\begin{equation}\label{eq: aux 9875}
\Prob\big\{\|G^{-k}\|_{HS}\geq t^k\,\sqrt{n}\big\}=O_k(1/t),\quad t\in[1,\infty).
\end{equation}

\bigskip

Next, we consider lower bounds for $\Prob\{\|G^{-k}\|\geq t^k\sqrt{n}\}$.
From now on we fix $t\geq 1$.
Let us start by recalling the formula for the joint distribution density of eigenvalues of $G G^\top$
(see, for example, \cite[formula~(4.5)]{EdelmanRao} or \cite[formula~(9)]{Chafai}):
\begin{equation}\label{eq: joint}
\rho(\lambda_1,\dots,\lambda_n):=
c(n)\,\exp\Big(-\frac{1}{2}\sum_{i=1}^n \lambda_i\Big)\,\prod\limits_{i<j}(\lambda_i-\lambda_j)\,\prod\limits_{i=1}^n \lambda_i^{-1/2},
\quad \lambda_1\geq\lambda_2\geq\dots\geq\lambda_n\geq 0,
\end{equation}
where $c(n)$ is a normalizing factor.
Denote
$$\mathcal L':=\Big\{(a_1,\dots,a_{n-1})\in\R^{n-1}:\;a_1\geq\dots\geq a_{n-1}\geq 0,\;\sum_{i=1}^{n-1} a_i^{-1}\leq 
4C_{\text{\tiny{\ref{thm: G_HS}}}}C t^2\,n\Big\},$$
where $C>0$ is the implicit constant from \eqref{eq:Szarek}, and $C_{\text{\tiny{\ref{thm: G_HS}}}}$
is taken from \eqref{thm: G_HS}.
Further, for any vector $a'=(a_1,\dots,a_{n-1})\in \mathcal L'$, let $u(a')\geq 0$ be the smallest non-negative integer such that
\begin{equation}\label{eq: aux 20985}
\rho\big(a_1,\dots,a_{n-1},4^{-u(a')-1}/(16C_{\text{\tiny{\ref{thm: G_HS}}}}^2C^4 t^2\,n)\big)
\leq 4\rho\big(a_1,\dots,a_{n-1},4^{-u(a')}/(16C_{\text{\tiny{\ref{thm: G_HS}}}}^2C^4 t^2\,n)\big).
\end{equation}
A simple analysis of formula \eqref{eq: joint} shows that $u(a')$ is well defined for any $a'\in\mathcal L'$.
Note that the definition of $u(a')$ implies that
\begin{equation}\label{eq: aux 0985209587}
\int\limits_{a_n=4^{-u(a')-1}/(16C_{\text{\tiny{\ref{thm: G_HS}}}}^2
C^4 t^2\,n)}^{4^{-u(a')}/(16C_{\text{\tiny{\ref{thm: G_HS}}}}^2 C^4 t^2\,n)}\rho(a_1,a_2,\dots,a_n)\,da_n
\geq \frac{1}{4}
\int\limits_{a_n=4^{-1}/(16C_{\text{\tiny{\ref{thm: G_HS}}}}^2
C^4 t^2\,n)}^{1/(16C_{\text{\tiny{\ref{thm: G_HS}}}}^2 C^4 t^2\,n)}\rho(a_1,a_2,\dots,a_n)\,d a_n.
\end{equation}
Now, we set
\begin{align*}
\mathcal L:=\big\{&(a_1,\dots,a_n)\in\R^n:\;a'=(a_1,\dots,a_{n-1})\in\mathcal L',\\
&a_n\in [4^{-u(a')-1}/(16C_{\text{\tiny{\ref{thm: G_HS}}}}^2C^4 t^2\,n),4^{-u(a')}/(16C_{\text{\tiny{\ref{thm: G_HS}}}}^2C^4 t^2\,n)]\big\}.
\end{align*}
It can be checked that $\mathcal L$ is a Borel set. Further, we clearly have
$$\int\limits_{(a_1,\dots,a_n)\in\mathcal L}\rho(a_1,\dots,a_n)\,d a_1\dots da_n=\Prob\big\{(s_1^2(G),\dots,s_{n}^2(G))
\in \mathcal L\big\}.$$
Let us show that the above quantity is bounded from below by $\widetilde c/t$
for a universal constant $\widetilde c>0$.
By combining~\eqref{eq:Szarek} with~\eqref{thm: G_HS}, we get that the event
$$
\big\{\|G^{-1}\|_{HS}^2\leq 16C_{\text{\tiny{\ref{thm: G_HS}}}}^2C^2
 t^2\,n\mbox{ and }\|G^{-1}\|\in [t\,n^{1/2},4C^2 t\,n^{1/2}]\big\}
$$
has probability at least $\frac{1}{Ct}-\frac{C}{4C^2 t}-\frac{C_{\text{\tiny{\ref{thm: G_HS}}}}}{4C_{\text{\tiny{\ref{thm: G_HS}}}}C t}
\geq \frac{1}{2Ct}$. On the other hand, everywhere on that event we have
$\sum_{i=1}^{n-1}s_i^{-2}(G)\leq 16C_{\text{\tiny{\ref{thm: G_HS}}}}^2C^2 t^2 n$
and $\frac{1}{16 C^4t^2n}\leq s_n^2(G)\leq \frac{1}{t^2n}$.
For any fixed $a_1\geq \dots\geq a_{n-1}\geq 0$, the density $\rho(a_1,a_2,\dots,a_{n-1},a_n)$,
viewed as a function of $a_n\in[0,a_{n-1}]$, is non-increasing. Hence, with $u(a')$ defined as above and in view of \eqref{eq: aux 0985209587}, we have
for every $a'=(a_1,\dots,a_{n-1})\in\mathcal L'$:
\begin{align*}
\int\limits_{a_n=4^{-u(a')-1}/(16C_{\text{\tiny{\ref{thm: G_HS}}}}^2
C^4 t^2\,n)}^{4^{-u(a')}/(16C_{\text{\tiny{\ref{thm: G_HS}}}}^2 C^4 t^2\,n)}\rho(a_1,a_2,\dots,a_n)\,da_n
\geq c
\int\limits_{a_n=1/(16 C^4 t^2n)}^{1/(t^2 n)}\rho(a_1,a_2,\dots,a_n)\,{\bf 1}_{\{a_{n}\leq a_{n-1}\}}\,da_n
\end{align*}
for some $c>0$, whence
$$
\int\limits_{(a_1,\dots,a_n)\in\mathcal L}\rho(a_1,\dots,a_n)\,d a_1\dots da_n\geq \frac{c/(2C)}{t}=:\frac{\widetilde c}{t}.
$$
Now, fix any $a'=(a_1,\dots,a_{n-1})\in\mathcal L'$, and apply the third assertion of Proposition~\ref{prop:HSbound}:
denoting by $T(a',s)$ the diagonal matrix with $T(a',s)_{jj}=a_j^{-1/2}$ (for $j<n$) and
the $(n,n)$--th entry equal to $s$, we get
\begin{align*}
\Big|\Big\{s:\;
&s^{-2}\in [4^{-u(a')-1}/(16C_{\text{\tiny{\ref{thm: G_HS}}}}^2C^4 t^2\,n),4^{-u(a')}/(16C_{\text{\tiny{\ref{thm: G_HS}}}}^2C^4 t^2\,n)],\\
&\E_W\brkb{\brkb{\brk{T(a',s) W}^{k-1} T(a',s)}_{nn}}^2\geq 
c_{k-1}\big(4^{u(a')}\cdot 16C_{\text{\tiny{\ref{thm: G_HS}}}}^2C^4 t^2\,n\big)^{k}n^{-k+1}\Big\}\Big|\\
&\hspace{3cm}\geq 2\cdot 2^{u(a')}C_{\text{\tiny{\ref{thm: G_HS}}}} C^2 t\,n^{1/2}.
\end{align*}
In view of \eqref{eq: aux 20985} and the lower bound for $\Prob\big\{(s_1^2(G),\dots,s_{n}^2(G))
\in \mathcal L\big\}$, the last inequality implies that
\begin{align*}
\Prob_{\Sigma}\Big\{&\|\Sigma^{-1}\|_{HS}^2\leq 4^{u(s_1^2(\Sigma),\dots,s_{n-1}^2(\Sigma))+1}
\cdot 16C_{\text{\tiny{\ref{thm: G_HS}}}}^2C^4 t^2\,n
+4C_{\text{\tiny{\ref{thm: G_HS}}}}C t^2\,n\quad\mbox{ and}\\
&\E_W\brkb{\brkb{\brk{\Sigma^{-1} W}^{k-1} \Sigma^{-1}}_{nn}}^2\geq 
c_{k-1}\big(4^{u(s_1^2(\Sigma),\dots,s_{n-1}^2(\Sigma))}\cdot 16C_{\text{\tiny{\ref{thm: G_HS}}}}^2C^4 t^2\big)^{k}n
\Big\}\geq \frac{c''}{t}
\end{align*}
for a universal constant $c''>0$.
The first assertion of Proposition~\ref{prop:HSbound} with $m=4$ and the last estimate yield
\begin{align*}
\Prob_{\Sigma}\Big\{&\E_W\brkb{\brkb{\brk{\Sigma^{-1} W}^{k-1} \Sigma^{-1}}_{nn}}^4
\leq C''_k \big(4^{u(s_1^2(\Sigma),\dots,s_{n-1}^2(\Sigma))}t^2\big)^{2k}n^2
\quad\mbox{ and}\\
&\E_W\brkb{\brkb{\brk{\Sigma^{-1} W}^{k-1} \Sigma^{-1}}_{nn}}^2\geq 
c_{k-1}\big(4^{u(s_1^2(\Sigma),\dots,s_{n-1}^2(\Sigma))}\cdot 16C_{\text{\tiny{\ref{thm: G_HS}}}}^2C^4 t^2\big)^{k}n
\Big\}\geq \frac{c''}{t}
\end{align*}
Applying the Paley--Zygmund inequality inside $\Prob_\Sigma\{\dots\}$, we get
$$
\Prob_{\Sigma}\Big\{\Prob_W\big\{\brkb{\brkb{\brk{\Sigma^{-1} W}^{k-1} \Sigma^{-1}}_{nn}}^2\geq 
c''_k t^{2k}n\big\}\geq \widetilde c_k
\Big\}\geq \frac{c''}{t}
$$
for some $c_k'',\widetilde c_k>0$ depending only on $k$,
whence
$$
\Prob\big\{\|G^{-k}\|^2\geq c''_k t^{2k}n\big\}\geq \frac{c'' \widetilde c_k}{t}.
$$
It remains to note that, together with the deterministic relation $\|G^{-k}\|_{HS}\geq \|G^{-k}\|$,
the above inequality and \eqref{eq: aux 9875} imply
$$
\Prob\big\{\|G^{-k}\|_{HS}\geq t^k\,\sqrt{n}\big\}=\Theta_k(1/t)\quad\mbox{ and }\quad
\Prob\big\{\|G^{-k}\|\geq t^{k}\sqrt{n}\big\}=\Theta_k(1/t),\quad t\in[1,\infty),
$$
and the theorem follows.
 
\subsection{Proof of Proposition~\ref{prop:HSbound}}

    Note that for any deterministic $n\times n$ matrix $B=(b_{ij})$ and any $k\in\N$,
    the $(i,j)$--th entry of $B^k$ can be expressed as 
    \begin{align} \label{eq: entryOfProducts}
        \brkb{B^k}_{ij} = 
        \sum_{\a \in [n]^{k-1}}
            b_{i\a_1}b_{\a_1\a_2}\cdots b_{\a_{k-1}j}.
    \end{align}
    Here, for $k=1$ we assume that $[n]^0$ consists of a single ``empty'' index vector $\alpha$.

    Let $P={\rm diag}(\delta_1,\dots,\delta_n)$ be a random matrix such that $\delta_i$
    are i.i.d.\ random signs $\pm 1$ jointly independent with $W$. Then $PW$ and $W$ 
    (hence, $\brk{TW}^{k}T$ and $\brk{TPW}^{k}T$) have
    the same distribution. 
    Applying \eqref{eq: entryOfProducts} to $T PW$ in place of $B$, we get for any $i,j\in[n]$: 
    \begin{align*}
        \brkb{\brk{T P W}^kT}_{ij} &= 
        \sum_{\a \in [n]^{k-1}} 
        \brk{\tau_i\tau_{\a_1}\cdots \tau_{\a_{k-1}}}
        \brk{\delta_i \delta_{\a_1}\cdots \delta_{\a_{k-1}}}
        w_{i\a_1}w_{\a_1\a_2}\cdots w_{\a_{k-1}j} \tau_{j}\\
        &= \tau_i\tau_j\delta_i\sum_{\a \in [n]^{k-1}} 
        \brk{\tau_{\a_1}\cdots \tau_{\a_{k-1}}}
        \brk{ \delta_{\a_1}\cdots \delta_{\a_{k-1}}}
        w_{i\a_1}w_{\a_1\a_2}\cdots w_{\a_{k-1}j},
    \end{align*}
    with the appropriate modification for the case $k=1$.
    To simplify the formulas, for any $m\geq 1$ and any index vector $\alpha \in [n]^{m(k-1)}$ we define
    $$
        w_{i,j,\alpha} := \prod\limits_{\ell=0}^{m-1} \big( w_{i\alpha_{\ell(k-1)+1}}
        \cdot w_{\alpha_{\ell(k-1)+1}\alpha_{\ell(k-1)+2}}\cdots w_{\alpha_{(\ell+1)(k-1)}j}\big).
    $$
    Then
    \begin{align*}
        \brkb{\brkb{\brk{T P W}^k T}_{ij}}^m = \tau_i^m\tau_j^m \sum_{\a\in [n]^{m(k-1)}}
        \brk{\prod_{\ell=1}^{m(k-1)}\tau_{\a_\ell}} \brk{\prod_{\ell=1}^{m(k-1)}\delta_{\a_\ell}} w_{i,j,\a},\quad m\geq 1.
    \end{align*}
    Note that for $m\geq 1$ and any given index vector
    $\alpha \in [n]^{m(k-1)}$, we have
    $\E_P \prod_{\ell=1}^{m(k-1)}\delta_{\a_\ell}=0$ if and only if there exists $h \in [n]$ such that $\abs{\set{\ell\,:\, \a_\ell = h}}$
    is odd. Let 
    $$ \Omega_m: =\set{ \a \in [n]^{m(k-1)}\,:\, 
    \forall h \in [n],\,\abs{\set{\ell\,:\, \a_\ell = h}} \text{ is even}}.
    $$
    Then
    \begin{align}\label{eq:aux 098709}
        \E_P \brk{\brkb{\brk{T P W}^k T}_{ij}}^m = \tau_i^m\tau_j^m \sum_{\a \in \Omega_m}
        \brk{\prod_{\ell=1}^{m(k-1)}\tau_{\a_\ell}}  w_{i,j,\a},\quad m\geq 1.
    \end{align}
    Next, observe that for any collection of $q$ non-negative random variables $X_1,\dots, X_{q}$ with 
    identical distributions we have
    \begin{align*}
        \E \prod_{\ell\in [q]}X_\ell \le \E\, \frac{1}{q!}\Big(\sum_{\ell \in [q]} X_\ell\Big)^{q}
        \leq \frac{1}{q!}\,\Big(\sum_{\ell \in [q]} \big(\E {X_\ell}^q\big)^{1/q}\Big)^{q}
        = \frac{q^q}{q!}\,{\E X_1}^{q},
    \end{align*}
    where we applied the triangle inequality for the $L_q$--norm in the second inequality. 
    Applying this relation to $w_{i,j,\a}$, we get
    \begin{align*}
        \abs{ \E w_{i,j,\a} } \le \frac{(mk)^{mk}}{(mk)!}\E |w_{11}|^{mk} \leq C_{k,m}n^{-km/2},\quad m\geq 1,\quad \a\in[n]^{m(k-1)},
    \end{align*}
    for some $C_{k,m}>0$ depending only on $k$ and $m$,
    where the last inequality follows by a standard moment estimates for one-dimensional projections of a vector uniformly distributed
    on $S^{n-1}$. 
    Combining the above estimates, we obtain
    \begin{align*}
        \E_W \E_P \brk{\brkb{\brk{T P W}^k T}_{ij}}^m &= 
        \E_W\,\tau_i^m\tau_j^m \sum_{\a \in \Omega_m}
        \brk{\prod_{\ell=1}^{m(k-1)}\tau_{\a_\ell}}  w_{i,j,\a}\\
        &\le C_{k,m}\tau_i^m\tau_j^m  \sum_{\a \in \Omega_m}
        \brk{\prod_{\ell=1}^{m(k-1)}\tau_{\a_\ell}} n^{-km/2}.
    \end{align*}
    Finally, for any even $m$ we construct a mapping $F_m$ from $\Omega_m$ to $[n]^{m(k-1)/2}$ as follows. Take any $\a\in \Omega_m$,
    and, at zeroth step, set $\gamma:=\a$.
    At step $1$, we set $\b_1 := \gamma_1$ and update the vector
    $\gamma$ by erasing both its first component and the component with the smallest index which is equal to $\b_1$.
    Thus, after the first step the vector $\gamma$ has length $m(k-1)-2$. At the second step,
    we set $\b_2:=\gamma_1$ and update
    $\gamma$ by erasing $\gamma_1$ and the first (other) component equal to $\b_2$.
    Thus, the length of $\gamma$ after the second step is $m(k-1)-4$.
    The validity of the procedure is guaranteed by the condition $\a\in \Omega_m$.
    After $m(k-1)/2$ steps we obtain a $m(k-1)/2$--dimensional vector $\b=(\b_1,\dots,\b_{m(k-1)/2})=:F(\a)$.
    It is not difficult to see that for every $\a \in \Omega_m$,
    $$
        \prod_{\ell\in [m(k-1)]} \tau_{\a_\ell}= \prod_{\ell\in [m(k-1)/2]}\tau_{F(\a)_\ell}^2.
    $$ 
    Therefore, for some $C_{k,m}'>0$ depending only on $k$ and $m$, we have
    \begin{align*}
        \E  \brk{\brkb{\brk{T W}^k T}_{ij}}^m&=
        \E_W \E_P \brk{\brkb{\brk{T P W}^k T}_{ij}}^m\\
        &\le C_{k,m}'\tau_i^m\tau_j^m  \sum_{\b \in [n]^{m(k-1)/2}}
        \brk{\prod_{\ell=1}^{m(k-1)/2}\tau_{\b_\ell}^2} n^{-km/2},
    \end{align*}
    giving the first assertion of the proposition.
    Letting $m=2$ and summing up over all $i\in[n]$ and $j\in[n]$, we obtain 
    $$
    \E \normb{\brk{T W}^k T}_{HS}^2 
    \le C_{k}'' n^{-k} \sum_{\b \in [n]^{k+1}} \prod_{i=1}^{k+1}{\tau_{\b_i}^2}
    = C_{k}'' n^{-k} \norm{T}_{HS}^{2k+2}
    $$
    for some $C_k''>0$ depending only on $k$, which gives the second assertion.

To prove the third assertion, we will use formula \eqref{eq:aux 098709},
which we will rewrite for $i=j$, $m=2$, and with the matrix $T$ replaced with $T(i,s)$.
We get
$$
\E \brk{\brkb{\brk{T(i,s) W}^k T(i,s)}_{ii}}^2 = s^4 \sum_{\a \in \Omega_m}
        \brk{\prod_{\ell=1}^{2(k-1)}(\tau_{\a_\ell}{\bf 1}_{\{\a_\ell\neq i\}}+s{\bf 1}_{\{\a_\ell= i\}})}  \E w_{i,i,\a}.
$$
Note that the above expression, viewed as a function of $s$, is a polynomial of degree $2k+2$,
and with the leading coefficient equal to $\E w_{ii}^{2k}=\Theta_k(n^{-k})$.
It follows immediately that on the interval $s\in[\tau_i/2,\tau_i]$, the polynomial is at least of order
$c_kn^{-k}\tau_i^{2k+2}$ on a set of Lebesgue measure $\tau_i/4$
(of course, we could write $(1-\delta)\tau_i/2$ for any constant $\delta>0$, at expense of decreasing $c_k$).
The result follows.

\subsection{Further remarks}
The corresponding problem for non-Gaussian matrices seems to be much more complicated due to the lack of rotational invariance.
It is natural to conjecture that for any $n\times n$ matrix $A$ with i.i.d.\ entries equidistributed with a random variable
$\xi$ of zero mean and unit variance,
$$
c_k\,t\leq \Prob\big\{s_{\min}(A^k)\leq t^k\,n^{-1/2}\big\}\leq C_k\,t,\quad 2e^{-c_k n}\leq t\leq 1,
$$
where $c_k,C_k>0$ may only depend on $k$ and the distribution of $\xi$ (and not on $n$).


\begin{thebibliography}{99}

\bibitem{Chafai}
{
D. Chafa\"i, Singular values of random matrices, 2009.
}

\bibitem{Edelman}
{
A. Edelman, Eigenvalues and condition numbers of random matrices, SIAM J. Matrix Anal. Appl. {\bf 9} (1988), no.~4, 543--560. MR0964668
}

\bibitem{EdelmanRao}
{
A. Edelman\ and\ N. R. Rao, Random matrix theory, Acta Numer. {\bf 14} (2005), 233--297. MR2168344
}

\bibitem{vNG}
{
J. von Neumann\ and\ H. H. Goldstine, Numerical inverting of matrices of high order, Bull. Amer. Math. Soc. {\bf 53} (1947), 1021--1099. MR0024235
}

\bibitem{Rud14}
{
M. Rudelson, Recent developments in non-asymptotic theory of random matrices, in {\it Modern aspects of random matrix theory},
83--120, Proc. Sympos. Appl. Math., 72, Amer. Math. Soc., Providence, RI. MR3288229
}

\bibitem{RV10}
{
M. Rudelson\ and\ R. Vershynin, Non-asymptotic theory of random matrices: extreme singular values,
in {\it Proceedings of the International Congress of Mathematicians. Volume III}, 1576--1602, Hindustan Book Agency, New Delhi. MR2827856
}

\bibitem{Smale}
{
S. Smale, On the efficiency of algorithms of analysis, Bull. Amer. Math. Soc. (N.S.) {\bf 13} (1985), no.~2, 87--121. MR0799791
}

\bibitem{Szarek}
{
S. J. Szarek, Condition numbers of random matrices, J. Complexity {\bf 7} (1991), no.~2, 131--149. MR1108773
}

\bibitem{Tikhomirov}
{
K. Tikhomirov, Invertibility via distance for non-centered random matrices with continuous distributions, arXiv:1707.09656.
}

\end{thebibliography}
\end{document}